\documentclass[10pt]{amsart}
\input{epsf}

\usepackage{amssymb,latexsym}
\topskip  \newtheorem{theorem}{Theorem}[section]

\newtheorem{corollary}[theorem]{Corollary}

\newtheorem{remark}[theorem]{Remark}
\newtheorem{lemma}[theorem]{Lemma}

\newtheorem{example}[theorem]{Example}

\begin{document}

\pagenumbering{arabic} \pagestyle{headings}

\def\LL{\frak L}
\def\AA{\frak A}
\def\BB{\frak B}
\def\HH{\frak H}
\def\TT{\frak T}

\def\pp{\overline{p}}
\def\qq{\overline{q}}
\def\zz{\overline{z}}
\def\NL{\mathcal{L}}
\def\sl{\sum\limits}
\def\xxj{\overline{x_j}}

\title{Average norms of polynomials}
\maketitle

\begin{center}Toufik Mansour \footnote{Research financed by EC's
IHRP Programme, within the Research Training Network "Algebraic
Combinatorics in Europe", grant HPRN-CT-2001-00272}
\end{center}
\begin{center}
Department of Mathematics, Chalmers University of Technology,
S-412~96 G\"oteborg, Sweden

{\tt toufik@math.chalmers.se}
\end{center}
\section*{Abstract}
In this paper we study the average $\NL_{2\alpha}$-norm over
$T$-polynomials, where $\alpha$ is a positive integer. More
precisely, we present an explicit formula for the average
$\NL_{2\alpha}$-norm over all the polynomials of degree exactly
$n$ with coefficients in $T$, where $T$ is a finite set of complex
numbers, $\alpha$ is a positive integer, and $n\geq0$. In
particular, we give a complete answer for the cases of Littlewood
polynomials and polynomials of a given height. As a consequence,
we derive all the previously known results for this kind of
problems, as well as many new results.

{\bf Keywords}: Littlewood polynomials, polynomials of height $h$,
$\NL_\alpha$-norm, generating function.
\section{Introduction}
The set of polynomials with special coefficients has given much
attention and there are a number of research old questions
concerning it. Erd\"os and Szekere; Hilbert; Littlewood; Prouhet,
Tarry, and Escott, are suggested a number of these questions
(see~\cite[Pages~5-7]{B2}). For example, Prouhet, Trray, and
Escott they asked to find a polynomial with integer coefficients
that is divisible by $(z-1)^n$ and has smallest sum of the
absolute values of the coefficients, and Erd\"os and Szekere they
asked to find the minimum of $||\prod_{j=1}^n(1-z^{a_j})||_\infty$
where the $a_j$ are positive integers, for given $n$. The problem
to find the maximum and the minimum norms of polynomials with
special coefficients is a specific old and difficult problem. In
the $\NL_4$-norm this problem is often called Golay's "Merit
Factor" problem. In the supremum norm this problem is due to
Littlewood. These problems are at least fifty years old and
neither solved. In this paper we interested in the following
problem (for particular cases see \cite[page~35]{B2}): Find the
average $\NL_{2\alpha}$-norm over polynomials of degree exactly
$n$ with coefficients in a given finite set $T$, where $\alpha$ is
a positive integer. In this paper, we find a complete answer for
this problem.

Let $T=\{x_1,x_2,\ldots,x_d\}$ be any finite set of complex
numbers. A polynomial $p(z)=a_nz^n+\cdots+a_1z+a_0$ is said to be
{\em $T$-polynomial} if $a_i\in T$ for all $i$, $0\leq i\leq n$.
We denote by $\TT_T(n)$ the set of all $T$-polynomials of degree
exactly $n$. For example, if $T=\{0,1,2\}$ then the set of
$T$-polynomials of degree exactly $1$ is given by $\TT_T(1)=\{z,
2z, z+1, 2z+1, z+2 ,2z+2\}$. The cardinality of this set,
$\TT_T(n)$, is denoted by $N_T(n)$. Clearly, for all $n\geq0$,
$$N_T(n)=\left\{\begin{array}{ll}
              d^{n+1},& 0\not\in T\mbox{ and }n\geq1,\\
              (d-1)d^n,& 0\in T\mbox{ and }n\geq1,\\
              d,& n=0.
              \end{array}\right.$$

A polynomial $p(z)$ is said to be {\em Littlewood polynomial} if
it $\{-1,1\}$-polynomial. A polynomial $p(z)$
is said to be {\em polynomial of height $h$} if it
$\{-h,-h+1,\ldots,h-1,h\}$-polynomial. For
example, $p(z)=z^2-z+1$ is a Littlewood polynomial of degree $2$,
and $p(z)=z^3-2z^2-1$ is a polynomial of height $2$ and degree
$3$.

Let $p(z)$ be any $T$-polynomial of degree exactly $n$. For any
positive integer $\alpha$, the {\em $\NL_\alpha$-norm} on the
boundary of the unit disk is defined by
$$||p||_{\alpha}=\left( \frac{1}{2\pi}\int_{0}^{2\pi}
|p(e^{i\theta})|^{\alpha}\,d\theta\right)^{\frac{1}{\alpha}}.$$
Let $f,g,h$ be any three real $\mathbb{R}$-polynomials, the {\em
$(f,g,h)$-average} over $T$-polynomials of degree exactly $n$ is
defined by
\begin{equation}\label{defe}
E_T(n;f,g,h)=\frac{1}{2\pi N_T(n)}\sum_{p\in\TT_T(n)}\int_0^{2\pi}
    h(e^{i\theta})f(p(e^{i\theta}))g(p(e^{-i\theta}))\,d\theta,
\end{equation}
for any $n\geq1$. We denote by $e_T(n,s,t,m)$ the
$(z^s,z^t,z^m)$-average over $T$-polynomials of degree exactly
$n$, where $m\in\mathbb{Z}$ and $s,t,n\geq0$.


We define the {\em average $\NL_\alpha$-norm} over
$T$-polynomials of degree exactly $n$ by
\begin{equation}\label{defa}
\mu^\alpha_{T}(n)=e_{T}(n;\alpha/2,\alpha/2,0)
=\frac{1}{N_T(n)}\sum_{p\in\TT_T(n)} ||p||_\alpha^\alpha
=\frac{1}{2\pi N_T(n)} \sum_{p\in\TT_n}
\int_0^{2\pi}|p(e^{i\theta})|^\alpha\, d\theta,
\end{equation}
for any positive integer $\alpha$. For
$\alpha=0$ we define $\mu_T^0(n)=1$ for all $n\geq1$.

Ones can asked to find an explicit formula for
$\mu^{2\alpha}_{T}(n)$, where $T$ is a finite set of complex numbers and $\alpha$ is a
positive integer (for particular cases, see the
research problem in \cite[Page~35]{B2}).
While the cases of Littlewood polynomials and polynomials of
height $h$ have attracted much attention (for example, see
\cite{B2,BC,NB}), the case of other sets $T$. The case of
$T=\{-1,1\}$, Littlewood polynomials, considered by several
authors as follows. In 1990, Newman and Byrnes~\cite{NB} found
$\mu^4_{\{-1,1\}}(n)=2n^2+3n+1$. In 2002, Borwein and
Choi~\cite{BC} they proved
$$\begin{array}{l}
\mu^6_{\{-1,1\}}(n)=6n^3+9n^2+4n+1,\\
\mu^8_{\{-1,1\}}(n)=24n^4+30n^3+4n^2+5n+4-3(-1)^n.
\end{array}$$
In the case $T=\{-1,0,1\}$, polynomials of height $1$,
Browein~\cite{Bt} proved $\mu^2_{\{-1,0,1\}}(n)=\frac{2}{3}(n+1)$,
$\mu^4_{\{-1,0,1\}}(n)=\frac{2}{9}(4n^2+7n+3)$, and
$\mu^6_{\{-1,0,1\}}(n)=\frac{2}{9}(8n^3+18n^2+13n+3)$. More
generally, Browein~\cite{Bt} found for any $h\geq 0$,
$$\begin{array}{l}
\mu^2_{\{-h,-h+1,\ldots,h-1,h\}}(n)=\frac{h(h+1)}{3}(n+1),\\
\mu^4_{\{-h,-h+1,\ldots,h-1,h\}}(n)=
\frac{h(h+1)}{45}(10h(h+1)n^2+(19h^2+19h-3)n+3(3h^2+3h-1)).
\end{array}$$

In this paper we suggest a general approach to study the average
$\NL_{2\alpha}$-norm over $T$-polynomials of degree exactly $n$,
for any positive integer $\alpha$ and any finite set $T$ of
complex numbers, which allows one to get an explicit expression
for $\mu^{2\alpha}_{T}(n)$. More precisely, we find an explicit
expression to the generating function for the sequence
$\{e_T(n;s,t,m)\}_{m,n,s,t}$. Using this generating function we
get an explicit expression for $e_T(n;s,t,m)$ in general, and
$\mu^{2\alpha}_T(n)=e_T(n,\alpha,\alpha,0)$ in particular, where
$m\in\mathbb{Z}$, $s,t,n\geq0$, $\alpha$ is a positive integer,
and $T$ is a finite set of complex numbers. As a consequence, we
derive all the previously known results for this kind of problems,
as well as many new results.

The main result of this paper can be formulated as follows. We
denote by $e_T(x,u,v,w)$ the generating function for the sequence
$\{e_T(n,s,t,m)\}_{n,s,t,m}$, that is,
$$e_T(x,u,v,w)=
\sl_{m\in\mathbb{Z}}\,\sl_{p,q,n\geq0}\,e_T(n,s,t,m)x^nu^sv^tw^m.$$

\begin{theorem}\label{maa}
Let $T=\{x_1,x_2,\ldots,x_d\}$ be a finite set of complex numbers.
Then the generating function $e_T(x,u,v,w)$ is given by
$$\sum_{n\geq1}\left[ \frac{1}{d^n}\sum_{j_1,\ldots,j_n=1}^d
\frac{x^{n-1}}{\left(1-x_{j_1}u-x_{j_2}uw^{-1}-\cdots-x_{j_n}uw^{-n+1}\right)
\left(1-\overline{x_{j_1}}v-\overline{x_{j_2}}vw-\cdots-\overline{x_{j_n}}vw^{n-1}\right)}\right].$$
Moreover, $e_T(n,s,t,m)$ is given by
$$\frac{1}{d^{n+1}}\sum_{j_1,\ldots,j_{n+1}=1}^d\sum_{s,t\geq0}
\sum_{\scriptsize\begin{array}{l}k_1+\cdots+k_{n+1}=s\\\ell_1+\cdots+\ell_{n+1}=t\\
\sum_{a=1}^{n+1}(a-1)(\ell_a-k_a)=m\end{array}}
\prod_{a=1}^{n+1}(x_{j_a}^{t_a}\overline{x_{j_a}}^{r_a})
\binom{s}{k_1,\ldots,k_{n+1}}\binom{t}{\ell_1,\ldots,\ell_{n+1}}.$$
\end{theorem}

Using Theorem~\ref{maa} we get an explicit expression for the
average $\NL_{2\alpha}$-norm over $T$-polynomials of degree
exactly $n$, namely $\mu_T^{2\alpha}$.

\begin{corollary}\label{ccmaa}\label{cc2}
Let $T=\{x_1,x_2,\ldots,x_d\}$ be a finite set of complex numbers.
Then the  average $\NL_{2\alpha}$-norm over $T$-polynomials of
degree exactly $n$, namely $\mu^{2\alpha}_T(n)$, is given by
$$\frac{1}{d^{n+1}}\sum_{j_1,\ldots,j_{n+1}=1}^d
\sum_{\scriptsize\begin{array}{c}k_1+\cdots+k_{n+1}=\alpha\\\ell_1+\cdots+\ell_{n+1}=\alpha\\
\sum_{a=1}^{n+1}(a-1)(\ell_a-k_a)=0\end{array}}
\prod_{a=1}^{n+1}(x_{j_a}^{k_a}\overline{x_{j_a}}^{\ell_a})
\binom{\alpha}{k_1,\ldots,k_{n+1}}\binom{\alpha}{\ell_1,\ldots,\ell_{n+1}}.$$
\end{corollary}
\begin{proof}
Applying Theorem~\ref{maa} and the identity
$e_T(n,\alpha,\alpha,0)=\mu^{2\alpha}_T(n)$ we get the desired result.
\end{proof}

The paper is organized as follows. The proof of our main result, Theorem~\ref{maa}, 
is presented in Section~\ref{sec2}. In
section~\ref{sec22} we present a general application for our
results. In particular, we give an explicit formulas up to
$\alpha=4$ where the sum of elements the set $T$ is $0$ (as in the
case of Littlewood polynomials and polynomials of height $h$). In
Section~\ref{sec3} we apply our main result for ceratin sets $T$, 
which allows us to get more details in the cases of Littlewood polynomials 
and polynomials of height $1$. Finally, in
Section~\ref{sec4} we suggest several directions to generalize the
results of the previous sections.
\section{Proofs}\label{sec2}
Let us start by introduce the quantity that plays the crucial role
in the proof of Theorem~\ref{maa}.

\begin{theorem}\label{th1}
Let $T=\{x_1,x_2,\ldots,x_d\}$ be a set of complex numbers,
$n\geq1$, $m$ any integer, and $s,t\geq0$. Then
\begin{equation}\label{eq1}
e_T(n,s,t,m)=\frac{1}{d}\sum_{j=1}^d\sl_{k=0}^a\sl_{\ell=0}^b
x_j^{s-k}\,\overline{x_j}^{\,t-\ell}\binom{s}{k}\binom{t}{\ell}e_T(n-1,k,\ell,m+k-\ell).
\end{equation}
\end{theorem}
\begin{proof}
Let $z=e^{i\theta}$ and $f_{s,t,m}(p(z))=z^mp^s(z)\pp^t(z)$, where $\pp(z)$ is the
conjugate polynomial of the polynomial $p(z)\in\TT_T(n)$. Since
for any $T$-polynomial $p(z)\in\TT_T(n)$ there exists an unique
polynomial $q\in\TT_T(n-1)$ and exists $j$, $1\leq j\leq d$, such
that $p(z)=zq(z)+x_j$, we have that
$$\sum_{p(z)\in\TT_T(n)}f_{s,t,m}(p(z))=z^m\sum_{j=1}^d\sum_{q(z)\in\TT_T(n-1)}
(zq(z)+x_j)^s(\zz\,\qq(z)+\overline{x_j})^t.$$ Using
$(x+y)^k=\sum_{j=0}^k \binom{k}{j}x^jy^{k-j}$ and $\zz=z^{-1}$
we get
$$
\sl_{p(z)\in\TT_T(n)}f_{s,t,m}(p(z))=\sl_{q(z)\in\TT_T(n-1)}\sum_{j=1}^d
\sl_{k=0}^s\sl_{\ell=0}^tx_j^{s-k}\,\overline{x_j}^{\,t-\ell}\binom{s}{k}
\binom{t}{\ell}f_{k,\ell,m+k-\ell}(q(z)).
$$
Therefore, using Definition~\ref{defe} we get the desired result.
\end{proof}

To present Recurrence~\ref{eq1} in terms of generating
functions we need the following lemma.

\begin{lemma}\label{lem1}
Let $F(x,y)=\sum_{s,t\geq0} d_{s,t}x^sy^t$ be a generating
function with two variables. Then
$$\sum_{s,t\geq0}x^sy^t\left(\sum_{k=0}^s\sum_{\ell=0}^t a^{s-k}b^{t-\ell}d_{k,\ell}
\binom{s}{k}\binom{t}{\ell}\right)=\frac{1}{(1-ax)(1-by)}F\left(\frac{x}{1-ax},\frac{y}{1-by}\right).$$
\end{lemma}
\begin{proof}
By definitions we have
$$\frac{1}{(1-ax)(1-by)}F\left(\frac{x}{1-ax},\frac{y}{1-by}\right)
=\sum_{s,t\geq0}a^{-s}b^{-t}d_{s,t}\frac{(ax)^s(by)^t}{(1-ax)^{s+1}(1-by)^{t+1}}.$$
Using $\frac{x^r}{(1-x)^r}=\sum_{n\geq0}\binom{n}{r}x^r$ we get
$$\frac{1}{(1-ax)(1-by)}F\left(\frac{x}{1-ax},\frac{y}{1-by}\right)=\sum_{s,t\geq0}
\left(\sum_{k,\ell\geq0}a^{k-s}b^{\ell-t}d_{s,t}\binom{k}{s}\binom{\ell}{t}x^ky^\ell\right),$$
equivalently,
$$\frac{1}{(1-ax)(1-by)}F\left(\frac{x}{1-ax},\frac{y}{1-by}\right)=
\sum_{s,t\geq0}x^sy^t\left(\sum_{k=0}^s\sum_{\ell=0}^t
a^{s-k}b^{t-\ell}d_{k,\ell} \binom{s}{k}\binom{t}{\ell}\right),$$
as claimed.
\end{proof}

\begin{remark}
Lemma~\ref{lem1} can be generalized as follows. Let
$$F(v_1,\ldots,v_k)=\sum_{s_1,s_2,\ldots,s_k\geq0}
e_{s_1,\ldots,s_k}\prod_{d=1}^kv_d^{s_d}$$ be a generating
function with $k$ variables. Then
$$\sum_{s_1,s_2,\ldots,s_k\geq0}\prod_{d=1}^kv_d^{s_d}
\left(\sum_{j_1=0}^{s_1}\sum_{{j_2}=0}^{s_2}\cdots\sum_{j_k}^{s_k}
e_{j_1,j_2,\ldots,j_k}
\prod_{d=1}^k\binom{s_d}{j_d}w_d^{s_d-j_d}\right)$$ is given by
$$\frac{1}{\prod_{d=1}^k(1-v_dw_d)}F\left(\frac{v_1}{1-w_1v_1},\frac{v_2}{1-w_2v_2},\ldots,\frac{v_k}{1-w_kv_k}\right).$$
\end{remark}

Now we are ready to prove our main result, namely
Theorem~\ref{maa}.

\begin{theorem}\label{thma}
Let $T=\{x_1,x_2,\ldots,x_d\}$ be a finite set of complex numbers.
Then the generating function $e_T(x,v,u,w)$ is given by
$$\sum_{n\geq1}\left[ \frac{1}{d^n}\sum_{j_1,\ldots,j_n=1}^d
\frac{x^{n-1}}{\left(1-x_{j_1}u-x_{j_2}uw^{-1}-\cdots-x_{j_n}uw^{-n+1}\right)
\left(1-\overline{x_{j_1}}v-\overline{x_{j_2}}vw-\cdots-\overline{x_{j_n}}vw^{n-1}\right)}\right].$$
\end{theorem}
\begin{proof}
If multiplying Equation~\ref{eq1} by $x^nu^sv^tw^m$, and summing
over all $m\in\mathbb{Z}$, $s,t\geq0$, and $n\geq1$, together with
using Lemma~\ref{lem1}, then we arrive to
$$\begin{array}{l}
e_T(x,u,v,w)-\sl_{m\in\mathbb{Z}}\sl_{s,t\geq0}e(0,s,t,m)u^sv^tw^m
=\dfrac{x}{d}\sl_{j=1}^d
\left[\frac{1}{(1-x_ju)(1-\overline{x_j}v)}\,e_T\left(x,
\frac{u}{w(1-x_ju)},\frac{wv}{1-\overline{x_j}v},w\right)\right].
\end{array}$$
On the other hand, by definitions we have that
$e_T(0,s,t,m)=\frac{1}{d}\sum_{j=1}^d\delta_mx_j^s\overline{x_j}^t$
for any $m\in\mathbb{Z}$ and $s,t\geq0$, where $\delta_m=1$ if
$m=0$, otherwise $\delta_m=0$. So
$$\sum_{m\in\mathbb{Z}}\sum_{s,t\geq0}e_T(0,s,t,m)u^sv^tw^m=
\frac{1}{d}\sum_{j=1}^d\frac{1}{(1-x_ju)(1-\overline{x_j}v)}.$$
Therefore, by combining the above two equations we get that
\begin{equation}\label{eqms}
e_T(x,u,v,w)=\frac{1}{d}\sum_{j=1}^d\frac{1}{(1-x_ju)(1-\overline{x_j}v)}
\left[1+xe_T\left(x,\frac{u}{w(1-x_ju)},\frac{wv}{1-x_jv},w\right)\right].
\end{equation}
An infinite number of applications of this identity concludes the
proof.
\end{proof}

\begin{remark}
Theorem~\ref{thma} yields the generating function $e_T(x,u,v,w)$
is symmetric under the translation $(u,v,w)\rightarrow
(v,u,w^{-1})$, that is, $e_T(x,u,v,w)=e_T(x,v,u,w^{-1})$.
\end{remark}

Theorem~\ref{thma} can be presented as follows.

\begin{corollary}\label{cc1}
Let $T=\{x_1,x_2,\ldots,x_d\}$ be a finite set of complex numbers.
Then the generating function $e_T(x,v,u,w)$ is given by
$$\sum_{\tiny n\geq1}
\sum_{\tiny j_1,\ldots,j_n=1}^d\sum_{\tiny s,t\geq0}
\sum_{\tiny\begin{array}{l}k_1+\cdots+k_n=s\\\ell_1+\cdots+\ell_n=t\end{array}}
\prod_{a=1}^n(x_{j_a}^{t_a}\overline{x_{j_a}}^{r_a})
\binom{s}{\tiny k_1,\ldots,k_n}\binom{t}{\tiny \ell_1,\ldots,\ell_n}
\frac{x^{n-1}u^sv^tw^{\left(\sl_{a=1}^n
(a-1)(r_a-t_a)\right)}}{d^n}.$$
\end{corollary}
\begin{proof}
Using Theorem~\ref{thma} we get that the generating function
$e_T(x,v,u,w)$ is given by
$$\sum_{n\geq1}\left[
\sum_{j_1,j_2,\ldots,j_n=1}^d\sum_{s,t\geq0}
\left(x_{j_1}+x_{j_2}w^{-1}+\cdots+x_{j_n}w^{-n+1}\right)^s
\left(\overline{x_{j_1}}+\overline{x_{j_2}}w+\cdots+\overline{x_{j_n}}w^{n-1}\right)^t\frac{x^{n-1}u^sv^t}{d^n}\right].$$
the rest is easy to check by the identity
$(a_1+\ldots+a_n)^s=\sum_{k_1+\cdots+k_n=s}\binom{s}{k_1,\ldots,k_n}\prod_{j=1}^na_j^{k_j}$.
\end{proof}

Let us denote by $\mu_T^{\alpha}(x)$ the generating function for
the sequence $\{\mu^{\alpha}_T(n)\}_{n\geq0}$, that is,
$\mu_T^\alpha(x)=\sum_{n\geq0}\mu_T^\alpha(n)x^n$.
Corollary~\ref{ccmaa} gives a complete answer to find the generating function 
$\mu_T^{2\alpha}(x)$  for any given a finite set $T$ and a positive integer $\alpha$.

\begin{example}
Using Corollary~\ref{cc2} for $\alpha=1$ we get that
    $$\mu_T^2(x)=\sum_{n\geq1}\sum_{j_1,\ldots,j_n=1}^d
    \left(\sum_{k=1}^n|x_{j_k}|^2\right)\frac{x^{n-1}}{d^n}=
    \sum_{n\geq1}\left(nd^{n-1}\sum_{j=1}^d|x_j|^2\right)\frac{x^{n-1}}{d^n},$$ so it
is easy to see that
$\mu_T^2(x)=\frac{1}{d(1-x)}\sum_{k=1}^n|x_{j_k}|^2+x\mu_T^2(x)$,
hence
    $$\mu_T^2(x)=\sum_{n\geq0} \left(\frac{n+1}{d}\sum_{j=1}^d|x_j|^2\right)x^n.$$
In particular, we have that $\mu_{\{-1,1\}}(n;2)=n+1$, and
$\mu_{\{-1,0,1\}}(n;2)=\frac{2}{3}(n+1)$.
\end{example}
Corollary~\ref{cc2} provide a finite algorithm for finding the
average $\mu_T^{2\alpha}(n)$ where $n$, $T$, and $\alpha$ are
given. This algorithm has been implemented in MAPLE, and yields
explicit results for given $n$, $T$, and $\alpha$ (see the tables
below).

\begin{remark}\label{rem1}
If looking at Corollary~\ref{cc2} carefully, then we understand
that the algorithm is work very slowly. This since we have to
consider $d^n$ possibilities for $j_i$ and
$\binom{n+\alpha-1}{\alpha}^2$ possibilities for $k_i$ and
$\ell_i$. Thus we have to consider
$d^n\binom{n+\alpha-1}{\alpha}^2$ of possibilites. Thus 
make our algorithm hard to obtain new results for $n$ large.
\end{remark}
\section{Exact formulas}\label{sec22}
Corollary~\ref{cc2} provide a close formula for finding the
average $\NL_{2\alpha}$-norm over $T$-polynomials of degree
exactly $n$ for any given $n\geq0$ and $\alpha\geq1$.
Remark~\ref{rem1} yields that the problem to find exact formula
for $\mu_T^{2\alpha}(n)$ with given only $\alpha\geq1$ it is a
hard problem by using Corollary~\ref{cc2}. Thus, we suggest here
another approach to find an explicit formula for
$\mu_T^{2\alpha}(n)$.

First let us denote by $e_T^{s,t}(x,w)$ the $(s+t)$-derivative of
the generating function $e_T(x,u,v,w)$ with respect $u^s$ and then
with respect $v^t$ at $u=v=0$, that is,
$$e_T^{s,t}=e_T^{s,t}(x,w)=\left.\frac{\partial^{s+t}}{\partial u^s\partial
v^t}\,e_T(x,u,v,w)\right|_{u=v=0}.$$ For any $s,t\geq0$, we define
    $$A_T^{s,t}=\sum_{j=1}^d x_j^s\xxj^t.$$
Now let us consider Equation~\ref{eqms}. This equation provide a
finite algorithm, {\em $\mu$-algorithm}, for finding
$e_T(n,s,t,m)$ in general, and $\mu_{T}^{2\alpha}(n)$ in
particular, since $s!\,t!\,e_T(n,s,t,m)$ is the coefficient of
$w^mx^n$ in the $(s+t)$-derivative of the generating function
$e_T(x,u,v,w)$ with respect to $u^s$ and then with respect $v^t$
at $u=v=0$, namely $e_T^{s,t}(x,w)$, and
$\mu_T^{2\alpha}(n)=e_T(n;\alpha,\alpha,0)$. Therefore, the
$\mu$-algorithm with input $\alpha$ and output $\mu_T^{2\alpha}(x)$
can be constructed as follows:

\begin{enumerate}
\item Apply the derivative operator with respect $u^s$ and then
with respect $v^t$ on Equation~\ref{eqms} for all $s,t$, where
$0\leq s,t\leq \alpha$.

\item Find explicitly $e_T^{s,t}$  for all $s,t$, where $0\leq
s,t\leq\alpha$. This by solving the system equations which is
obtained from step 1.

\item Find $\mu_T^{2\alpha}(x)$, which is the free coefficient of
$w$ in $e_T^{\alpha,\alpha}(x,w)$.
\end{enumerate}
This algorithm has been implemented in MAPLE
, and yields explicit results 
for given $\alpha$. Below we present several explicit calculations.

\subsection{Formula for $\mu^2_T(n)$} Let us start by apply the
$\mu$-algorithm for $\alpha=1$. The first step of the
$\mu$-algorithm gives
$$\left\{\begin{array}{l}
e_T^{0,0}=1+xe_T^{0,0},\\
\\
e_T^{0,1}=\frac{1}{d}A_T^{0,1}+\frac{x}{d}A_T^{0,1}e_T^{0,0}+\frac{x}{w}e_T^{0,1},\\
\\
e_T^{1,0}=\frac{1}{d}A_T^{1,0}+\frac{x}{d}A_T^{1,0}e_T^{0,0}+\frac{x}{w}e_T^{1,0},\\
\\
e_T^{1,1}=\frac{1}{d}A_T^{1,1}+\frac{x}{d}A_T^{1,1}e_T^{0,0}+\frac{xw}{d}A_T^{1,0}e_T^{0,1}+
\frac{x}{dw}A_T^{0,1}e_T^{1,0}+xe_T^{1,1}.
\end{array}\right.$$
Equivalently (the second step of the $\mu$-algorithm),
$$e_T^{0,0}=\frac{1}{1-x},\quad
e_T^{1,0}=\frac{A_T^{1,0}}{d(1-x)(1-xw^{-1})},\quad
e_T^{0,1}=\frac{A_t^{0,1}}{d(1-x)(1-xw)},$$
and
$$\begin{array}{l}
e_T^{1,1}= \frac{1}{d(1-x)^2}A_T^{1,1}
+\left(\frac{xw}{d^2(1-x)^2(1-xw)}
+\frac{xw^{-1}}{d^2(1-x)^2(1-xw^{-1})}\right)A_T^{1,0}A_T^{0,1}.
\end{array}$$
Therefore, the third step of the $\mu$-algorithm gives
$\mu_T^2(x)$, which is the free coefficient of $w$ in
$e_T^{1,1}(x,w)$. Hence, we get the following result.

\begin{corollary}\label{case2}
We have
$$\mu_T^2(x)=\frac{1}{d(1-x)^2}A_T^{1,1}=\frac{1}{d(1-x)^2}\sum_{j=1}^d|x_j|^2.$$
Moreover, for all $n\geq0$,
$$\mu_T^2(n)=\frac{n+1}{d}A_T^{1,1}=\frac{n+1}{d}\sum_{j=1}^d|x_j|^2.$$
\end{corollary}

For example, in the case of Littlewood polynomials, namely
$T=\{-1,1\}$, we have that $\mu_T^2(n)=n+1$, and in the case of
polynomials of height $h$, namely, $T=\{-h,-h+1,\ldots,h-1,h\}$,
we have that $\mu_T^2(n)=\frac{h(h+1)}{3}(n+1)$.

\subsection{Formula for $\mu^4_T(n)$} Again, using the
$\mu$-algorithm for $\alpha=2$ we get that 

$\begin{array}{rl}
e_T^{0,2}&=\frac{2}{d(1-x)(1-xw^2)}A_T^{0,2}+\frac{2xw}{d(1-xw^2)}e_T^{0,1}A_T^{0,1},\\
\end{array}$

$\begin{array}{rl}
e_T^{2,0}&=\frac{2}{d(1-x)(1-xw^{-2})}A_T^{2,0}+\frac{2x}{dw(1-xw^{-2})}e_T^{1,0}A_T^{1,0},\\
\end{array}$

$\begin{array}{rl}
e_T^{1,2}&=\frac{2}{d(1-x)(1-xw)}A_T^{1,2}
+\frac{4xw}{d(1-xw)}e_T^{0,1}A_T^{1,1}
+\frac{xw^2}{d(1-xw)}e_T^{0,2}A_T^{1,0}
+\frac{2x}{dw(1-xw)}e_T^{1,0}A_T^{0,2}\\
&+\frac{4x}{d(1-xw)}e_T^{1,1}A_T^{0,1},\\
\end{array}$

$\begin{array}{rl}
e_T^{2,1}&=\frac{2}{d(1-x)(1-xw^{-1})}A_T^{2,0}
+\frac{4x}{dw(1-xw^{-1})}e_T^{1,0}(x,w)A_T^{1,1}
+\frac{x}{dw^2(1-xw^{-1})}e_T^{0,2}A_T^{0,1}
+\frac{2xw}{d(1-xw^{-1})}e_T^{1,0}A_T^{2,0}\\
&+\frac{4x}{d(1-xw^{-1})}e_T^{1,1}A_T^{1,0},\\
\end{array}$

$\begin{array}{rl}
e_T^{2,2}&=2\frac{4}{d(1-x)^2}A_T^{2,2}
+\frac{8xw}{d(1-x)}e_T^{0,1}A_T^{2,1}
+\frac{2xw^2}{d(1-x)}e_T^{0,2}A_T^{2,0}
+\frac{2x}{dw^2(1-x)}e_T^{2,0}A_T^{0,2}
+\frac{8x}{dw(1-x)}e_T^{1,0}A_T^{1,2}\\
&+\frac{16x}{d(1-x)}e_T^{1,1}A_T^{1,1}
+\frac{4xw}{d(1-x)}e_T^{1,2}A_T^{1,0}
+\frac{4x}{dw(1-x)}e_T^{2,1}A_T^{0,1}
\end{array}$

Solving this linear system in $e_T^{s,t}$ where $0\leq s,t\leq 2$, we get
an explicit expression for $e_T^{2,2}(x,w)$ (it is long to present
here). Using this expression we find the free coefficient $w$ of
$e_T^{2,2}(x,w)$, hence we get the following result.
\begin{corollary}\label{case4}
We have
$$\begin{array}{ll}
\mu_T^4(x)&=\frac{1}{d(1-x)^2}A_T^{2,2}+
\frac{4x}{d^2(1-x)^3}\left(A_T^{1,1}\right)^2
+\frac{2x^2(1+x)^2}{d^3(1-x^2)^3}\left[\left(A_T^{1,0}\right)^2A_T^{0,1}+\left(A_T^{0,1}\right)^2A_T^{1,0}\right]\\
&+\frac{8x^3}{d^4(1-x)^4(1+x)}\left(A_T^{1,0}\right)^2\left(A_T^{0,1}\right)^2.
\end{array}$$ Moreover, for all $n\geq0$,
$$\begin{array}{ll}
\mu^4_T(n)&=\frac{2\left(A_T^{1,0}\right)^2\left(A_T^{0,1}\right)^2}{3d^4}n^3
+\left\{\frac{2\left(A_T^{1,1}\right)^2}{d^2}+
\frac{\left[\left(A_T^{1,0}\right)^2A_T^{0,1}+\left(A_T^{0,1}\right)^2A_T^{1,0}\right]}{2d^3}
-\frac{\left(A_T^{1,0}\right)^2\left(A_T^{0,1}\right)^2}{d^4}\right\}n^2\\[2pt]
&+\left\{\frac{2\left(A_T^{1,1}\right)^2}{d^2}+\frac{A_T^{2,2}}{d}-
\frac{2\left(A_T^{1,0}\right)^2\left(A_T^{0,1}\right)^2}{3d^4}\right\}n
+\frac{A_T^{2,2}}{d}+
\frac{(1-(-1)^n)\left(A_T^{1,0}\right)^2\left(A_T^{0,1}\right)^2}{2d^4}\\[2pt]
&-\frac{(1-(-1)^n)\left[\left(A_T^{1,0}\right)^2A_T^{0,1}+\left(A_T^{0,1}\right)^2A_T^{1,0}\right]}{3d^3}.
\end{array}$$
\end{corollary}

For example, in the case of Littlewood polynomials we have that
$\mu^4_{\{-1,1\}}(n)=2n^2+3n+1$, and in the case of polynomials of
height $h$ we get that
$$\mu^4_{\{-h,-h+1,\ldots,h-1,h\}}(n;4)=\frac{h(h+1)}{45}(10h(h+1)n^2+(19h^2+19h-3)n+3(3h^2+3h-1)).$$

\subsection{Formula for $\mu_T^{2\alpha}(n)$ where $A_T^{1,0}=A_T^{0,1}=0$}
Similarly to the previous subsection, our results can be extended
to the case of $\mu_T^6(n)$. Since the answers become very
cumbersome, we present here only the simplest case when
$A_T^{1,0}=\sum_{j=1}^dx_j=0$. Therefore, if we apply our approach
for finding $\mu_T^{2\alpha}(n)$ where $T=\{x_1,x_2,\ldots,x_d\}$
such that $A_T^{1,0}=A_T^{0,1}=0$ we get the following results.

\begin{corollary}\label{case00}
Let $T=\{x_1,x_2,\ldots,x_d\}$ such that $\sum_{j=1}^dx_j=0$. Then

{\rm(i)} 

$\mu_T^2(x)=\frac{1}{d(1-x)^2}A_T^{1,1}$.

{\rm(ii)}

$\mu_T^4(x)=\frac{1}{d(1-x)^2}A_T^{2,2}+\frac{4x}{d^2(1-x)^3}\left(A_T^{1,1}\right)^2$.

{\rm(iii)}

$\mu_T^6(x)=\frac{1}{d(1-x)^2}A_T^{3,3}+\frac{18x}{d^2(1-x)^3}A_T^{1,1}A_T^{2,2}+\frac{36x^2}{d^3(1-x)^4}(A_T^{1,1})^3$.

{\rm(iv)} 

$\begin{array}{ll}
\mu_t^8(x)&=\frac{1}{d(1-x)^2}A_T^{4,4}+\frac{32x}{d^2(1-x)^3}A_T^{1,1}A_T^{3,3}
+\frac{36x^2}{d^2(1-x)^4}(A_T^{2,2})^2+\frac{432x^2}{d^3(1-x)^4}A_T^{2,2}(A_T^{1,1})^2\\
&+\frac{72x^4(3-2x-2x^2+3x^3-x^4)}{d^4(1-x)^4(1+x)}(A_T^{0,2})^2(A_T^{2,0})^2
+\frac{576x^3}{d^4(1-x)^5}(A_T^{1,1})^4\\
&+\frac{48x^3}{d^3(1-x)^2(1-x^3)}(A_T^{0,3}A_T^{2,0}A_T^{2,1}+A_T^{3,0}A_T^{0,2}A_T^{1,2})\\
&+\frac{72x^2}{d^3(1-x)^2(1-x^2)}((A_T^{1,2})^2A_T^{2,0}+(A_T^{2,1})^2A_T^{0,2})\\
&+\frac{6x^2}{d^3(1-x)^2(1-x^2)}((A_T^{2,0})^2A_T^{0,4}+(A_T^{0,2})^2A_T^{4,0}).
\end{array}$
\end{corollary}

We have two remarks. The first one is the algorithm more easy to
run under the condition of $A_T^{1,0}=A_T^{0,1}=0$. This, since
the corresponding generating function become more simple. For example, 
see Corollary~\ref{case4} and
(ii) in  Corollary~\ref{case00}. The second one is the generating function
$\mu_T^{2\alpha}(x)$ is symmetric on $x_1,x_2,\ldots,x_d$. This
since, we obtained the formula for $\mu_T^{2\alpha}(x)$ from
Equation~\ref{eqms} which is also symmetric on $x_1,x_2,\ldots,x_d$.
\section{Applications}\label{sec3}
In this section we discus a particular cases of Theorem~\ref{maa}
and Corollary~\ref{ccmaa}. More precisely, we apply
Corollary~\ref{ccmaa} for ceratin finite sets of complex numbers
(special real numbers), and then we present an explicit formulas
for $\mu_T^{2\alpha}$ for certain positive integers $\alpha$.

\subsection{Littlewood polynomials} Here we suggest formulas for
$\mu_T^{2\alpha}(n)$, where $T=\{-1,1\}$, namely the case of
Littlewood polynomials. Corollary~\ref{ccmaa} for $T=\{-1,1\}$
gives the following result.
\begin{corollary}\label{ccxa}
The  average $\NL_{2\alpha}$-norm over Littlewood polynomials of
degree exactly $n$, namely $\mu^{2\alpha}_{\{-1,1\}}(n)$, is given
by
$$\frac{1}{2^{n+1}}\sum_{j_1,\ldots,j_{n+1}=1}^2
\sum_{\scriptsize\begin{array}{c}k_1+\cdots+k_{n+1}=\alpha\\\ell_1+\cdots+\ell_{n+1}=\alpha\\
\sum_{s=1}^{n+1}(s-1)(\ell_s-k_s)=0\end{array}}(-1)^{\sum_{s=1}^{n+1}j_s(k_s+\ell_s)}
\binom{\alpha}{k_1,\ldots,k_{n+1}}\binom{\alpha}{\ell_1,\ldots,\ell_{n+1}}.$$
\end{corollary}
Using Lemma \ref{ccxa} we quickly generate the numbers
$\mu_{\{-1,1\}}^\alpha(n)$; the first few of these numbers are
given in Table \ref{tbx1}.
\begin{table}[ht]
  $$
  \begin{array}{rrrrrrrr}
    n\backslash \alpha & 0 & 1 & 2 & 3  & 4  & 5  \\
    \hline
    0\;\;  & 1    & 1    & 1    &  1   &   1  &   1  \\
    1\;\;  & 1    & 2    & 6    & 20   &  70  &  252 \\
    2\;\;  & 1    & 3    &15    & 93   & 651  & 4913 \\
    3\;\;  & 1    & 4    &28    &256   &2812  & 35024  \\
    4\;\;  & 1    & 5    &45    &545   &8149  &143945 \\
    5\;\;  & 1    & 6    &66    &996   &18882 &433116 \\
    6\;\;  & 1    & 7    &81    &1645  &37759 &1062697 \\
    7\;\;  & 1    & 8    &120   &2528  &68152 &2272128\\
    8\;\;  & 1    & 9    &153   &3681  &113961&4385969\\
    9\;\;  & 1    &10    &190   &5140  &179710&7839260\\
   10\;\;  & 1    &11    &231   &6941  &270451&13178561\\
   \hline
  \end{array}
  $$
  \caption{Values of $\mu_{\{-1,1\}}^{2\alpha}(n)$.}
  \label{tbx1}
\end{table}

On the other hand, if applying Equation~\ref{eqms} for
$T=\{-1,1\}$ we get that
\begin{equation}\label{eqx1}
\begin{array}{ll}
e_{\{-1,1\}}(x,u,v,w)&=\frac{1}{2(1-u)(1-v)}
\left[1+xe_{\{-1,1\}}\left(x,\frac{u}{w(1-u)},\frac{wv}{1-v},w\right)\right]\\
&+\frac{1}{2(1+u)(1+v)}
\left[1+xe_{\{-1,1\}}\left(x,\frac{u}{w(1+u)},\frac{wv}{1+v},w\right)\right].
\end{array}
\end{equation}
Now, if applying $\mu$-algorithm on Equation~\ref{eqx1} 
for $T=\{-1,1\}$ and $\alpha=0,1,\ldots,5$, 
then we get the following result.

\begin{corollary}
For all $n\geq0$,

{\rm(i)} $\mu^0_{\{-1,1\}}(n)=1$,

{\rm(ii)} $\mu^2_{\{-1,1\}}(n)=n+1$,

{\rm(iii)} $\mu^4_{\{-1,1\}}(n)=2n^2+3n+1$,

{\rm(iv)} $\mu^6_{\{-1,1\}}(n)=6n^3+9n^2+4n+1$,

{\rm(v)} $\mu^8_{\{-1,1\}}(n)=24n^4+30n^3+4n^2+5n+4-3(-1)^n$,

{\rm(vi)}
$\mu^{10}_{\{-1,1\}}(n)=120n^5+150n^4-350n^3+265n^2+281n-144-5(-1)^n(15n-29)$.
\end{corollary}

\subsection{Polynomials of height $1$} Here we suggest, in the case of polynomials of height
$1$, some formulas for $\mu_{\{-1,0,1\}}^{2\alpha}(n)$. Using
Lemma \ref{ccxa} we quickly generate the numbers
$\mu_{\{-1,0,1\}}^\alpha(n)$; the first few of these numbers are
given in Table \ref{tbx2}.
\begin{table}[ht]
  $$
  \begin{array}{rrrrrrrr}
    n\backslash \alpha & 0 & 1 & 2 & 3  & 4  & 5  \\
    \hline
    &&&&&&\\[-7pt]
    0\;\;  & 1    & \frac{2}{3}    & \frac{6}{9}    &  \frac{6}{9}   &   \frac{6}{9}  &   \frac{18}{27}
    \\[3pt]
    1\;\;  & 1    & \frac{4}{3}    & \frac{28}{9}   &  \frac{84}{9}  &   \frac{284}{9}& \frac{3036}{27}
    \\[3pt]
    2\;\;  & 1    & \frac{6}{3}    & \frac{66}{9}   &  \frac{330}{9} &  \frac{2018}{9}& \frac{42334}{27}
    \\[3pt]
    3\;\;  & 1    & \frac{8}{3}    & \frac{120}{9}  & \frac{840}{9}  &  \frac{7480}{9}& \frac{239832}{27}
    \\[3pt]
    4\;\;  & 1    & \frac{10}{3}   & \frac{190}{9}  & \frac{1710}{9} & \frac{19902}{9}& \frac{856010}{27}
    \\[3pt]
    5\;\;  & 1    & \frac{12}{3}   &  \frac{276}{9} & \frac{3036}{9} & \frac{43604}{9}& \frac{2348788}{27}
    \\[3pt]
    6\;\;  & 1    & \frac{14}{3}   & \frac{378}{9}  & \frac{4914}{9} & \frac{83866}{9}& \frac{5410646}{27}
    \\[3pt]
    7\;\;  & 1    & \frac{16}{3}   & \frac{496}{9}  & \frac{7440}{9} &\frac{147056}{9}& \frac{11040304}{27}
    \\[3pt]
    8\;\;  & 1    & \frac{18}{3}   & \frac{630}{9}  & \frac{10710}{9}&\frac{240502}{9}& \frac{20567042}{27}
    \\[3pt]
    9\;\;  & 1    & \frac{20}{3}   & \frac{780}{9}  & \frac{14820}{9}&\frac{372620}{9}& \frac{35735180}{27}
    \\[3pt]
    10\;\; & 1    & \frac{22}{3}   & \frac{946}{9}  & \frac{19866}{9}&\frac{552786}{9}& \frac{58715598}{27}
    \\[3pt]
    \hline
  \end{array}
  $$
  \caption{Values of $\mu_{\{-1,0, 1\}}^{2\alpha}(n)$.}
  \label{tbx2}
\end{table}

Applying Equation~\ref{eqms} for $T=\{-1,0,1\}$ then we get that
\begin{equation}\label{eqx2}
\begin{array}{ll}
e_{\{-1,0,1\}}(x,u,v,w)&=\frac{1}{3}\left[1+xe_{\{-1,0,1\}}\left(
x,\frac{u}{w},wv,w\right)\right]\\
&+\frac{1}{3(1-u)(1-v)}
\left[1+xe_{\{-1,0,1\}}\left(x,\frac{u}{w(1-u)},\frac{wv}{1-v},w\right)\right]\\
&+\frac{1}{3(1+u)(1+v)}
\left[1+xe_{\{-1,0,1\}}\left(x,\frac{u}{w(1+u)},\frac{wv}{1+v},w\right)\right].
\end{array}
\end{equation}
Now, by applying $\mu$-algorithm on Equation~\ref{eqx2} for $T=\{-1,0,1\}$ and
$\alpha=0,1,\ldots,5$ we get as follows.

\begin{corollary}
For all $n\geq0$,

{\rm(i)} $\mu^0_{\{-1,0,1\}}(n)=1$,

{\rm(ii)} $\mu^2_{\{-1,0,1\}}(n)=\frac{2}{3}(n+1)$,

{\rm(iii)} $\mu^4_{n,\{-1,0,1\}}(n)=\frac{2}{9}(4n^2+7n+3)$,

{\rm(iv)} $\mu^6_{\{-1,0,1\}}(n)=\frac{2}{9}(8n^3+18n^2+13n+3)$,

{\rm(v)}
$\mu^8_{\{-1,0,1\}}(n)=\frac{2}{27}(64n^4+176n^3+128n^2+37n+15-6(-1)^n)$,

{\rm(vi)}
$\mu^{10}_{\{-1,0,1\}}(n)=\frac{2}{81}(640n^5+2400n^4+630n^2+1337n-363-30(-1)^n(10n-13)$.
\end{corollary}

\subsection{Polynomials with coefficients $0,1$}
Here we suggest another case, the case of polynomials with
coefficients $0,1$, which allows us to find explicit formulas for
$\mu_{\{0,1\}}^{2\alpha}(n)$. Using Lemma \ref{ccxa} we quickly
generate the numbers $\mu_{\{0,1\}}^\alpha(n)$; the first few of
these numbers are given in Table \ref{tbx3}.

\begin{table}[ht]
  $$
  \begin{array}{rrrrrrrr}
    n\backslash \alpha & 0 & 1 & 2 & 3  &   \\
    \hline
    &&&&&\\[-7pt]
    0\;\;  & 1   & \frac{1}{2}    & \frac{1}{2}    &  \frac{4}{8}   &
    \\[3pt]
    1\;\;  & 1   & \frac{2}{2}    & \frac{4}{2}    & \frac{44}{8}   &
    \\[3pt]
    2\;\;  & 1   & \frac{3}{2}    & \frac{10}{2}   & \frac{204}{8}  &
    \\[3pt]
    3\;\;  & 1   & \frac{4}{2}    & \frac{19}{2}   & \frac{592}{8}  &
    \\[3pt]
    4\;\;  & 1   & \frac{5}{2}    & \frac{32}{2}   & \frac{1397}{8} &
    \\[3pt]
    5\;\;  & 1   & \frac{6}{2}    & \frac{49}{2}   & \frac{2826}{8} &
    \\[3pt]
    6\;\;  & 1   & \frac{7}{2}    & \frac{71}{2}   & \frac{5206}{8} &
    \\[3pt]
    7\;\;  & 1   & \frac{8}{2}    & \frac{98}{2}   & \frac{8876}{8} &
    \\[3pt]
    8\;\;  & 1   & \frac{9}{2}    & \frac{131}{2}  & \frac{14334}{8}&
    \\[3pt]
    9\;\;  & 1   & \frac{10}{2}   & \frac{170}{2}  & \frac{22084}{8}&
    \\[3pt]
    10\;\; & 1   & \frac{11}{2}   & \frac{216}{2}  & \frac{32828}{8}&
    \\[3pt]
    \hline
  \end{array}
  $$
  \caption{Values of $\mu_{\{0,1\}}^{2\alpha}(n)$.}
  \label{tbx3}
\end{table}

Applying Equation~\ref{eqms} for $T=\{0,1\}$ then we get that
\begin{equation}\label{eqx3}
\begin{array}{ll}
e_{\{0,1\}}(x,u,v,w)&=\frac{1}{2}\left[1+xe_{\{0,1\}}\left(
x,\frac{u}{w},wv,w\right)\right]+\frac{1}{2(1-u)(1-v)}
\left[1+xe_{\{0,1\}}\left(x,\frac{u}{w(1-u)},\frac{wv}{1-v},w\right)\right].
\end{array}
\end{equation}
Now, if applying $\mu$-algorithm on Equation~\ref{eqx3} for $T=\{0,1\}$ and
$\alpha=0,1,\ldots,5$ we get the following result.

\begin{corollary}
For all $n\geq0$,

{\rm(i)} $\mu^0_{\{0,1\}}(n)=1$,

{\rm(ii)} $\mu^2_{\{0,1\}}(n)=\frac{1}{2}(n+1)$,

{\rm(iii)}
$\mu^4_{\{0,1\}}(n)=\frac{1}{96}(4n^3+54n^2+92n+45+3(-1)^n)$,

{\rm(iv)}
$\mu^6_{\{0,1\}}(n)=\frac{1}{2560}(22n^5+460n^4+3100n^3+5600n^2+4143n+1130+75(-1)^n(3n+2))$.
\end{corollary}

\subsection{Example of polynomials with coefficients $0,i$ where $i^2=-1$}
Here we extend the examples which presented by several authors
(see the first section) to the case of sets of complex numbers.
For example if we apply the main results of the pervious sections
for $T=\{0,i\}$ we get the following result.

\begin{corollary}
For all $n\geq0$,

{\rm(i)} $\mu^0_{\{0,i\}}(n)=1$,

{\rm(ii)} $\mu^2_{\{0,i\}}(n)=\frac{1}{2}n+\frac{1}{2}$,

{\rm(iii)}
$\mu^4_{\{0,i\}}(n)=\frac{1}{96}(4n^3+54n^2+92n+45+3(-1)^n)$,

{\rm(iv)} $\begin{array}{ll}
\mu^6_{\{0,i\}}(n)&=\dfrac{1}{5120}(46n^5+890n^4+6320n^3+11200n^2+7789n+2455\\
&\qquad\qquad+15(-1)^n(23+29n)-60i^n(2+3i-in)+60(-i)^n(-2+3i+in)),\end{array}$

where $i^2=-1$.
\end{corollary}
\section{Further results}\label{sec4}
In this section we suggest several directions to generalize the results of the
previous sections.

\subsection{Average norms of $T$-polynomials with weights}
The first of these directions is to obtain an exact formula for
the average $\NL_{2\alpha}$-norm over $T$-polynomials of degree
exactly $n$ with weight $z^m$ for given $\alpha$, $n$, and $m$.
Let us define,
\begin{equation}\label{defga}
\mu^\alpha_{T}(n;m)=e_{T}(n;\alpha/2,\alpha/2,m)=\frac{1}{2\pi
N_T(n)} \sum_{p\in\TT_n}
\int_0^{2\pi}e^{im\theta}|p(e^{i\theta})|^\alpha\, d\theta,
\end{equation}
Clearly, $\mu^\alpha_T(n)=\mu^\alpha_T(n;0)$ for all $n$, and
$\alpha>0$. Now, ones can asked the following general problem.
Find $\mu^{2\alpha}_T(n;m)$ for any $\alpha>0$, $n\geq0$, and
$m\in\mathbb{Z}$. We can give a complete answer for this problem
by using Theorem~\ref{maa} and the definitions. For example, the
following result is true.

\begin{theorem}\label{gen1}
Let $T=\{x_1,x_2,\ldots,x_d\}$ be any finite set of complex
numbers. The generating function
$\sum_{n\geq0}\sum_{m\in\mathbb{Z}}\mu^{2\alpha}_T(n;m)x^nw^m$ is
given by
$$\sum_{n\geq1}
\sum_{j_1,\ldots,j_n=1}^d
\sum_{\scriptsize\begin{array}{l}k_1+\cdots+k_n=\alpha\\
\ell_1+\cdots+\ell_n=\alpha\end{array}}
\prod_{a=1}^n(x_{j_a}^{t_a}\overline{x_{j_a}}^{r_a})
\binom{\alpha}{k_1,\ldots,k_n}\binom{\alpha}{\ell_1,\ldots,\ell_n}
\frac{x^{n-1}w^{\left(\sl_{a=1}^n
(a-1)(r_a-t_a)\right)}}{d^n}.$$ Moreover, $\mu^{2\alpha}_T(n;m)$
is given by
$$\frac{1}{d^{n+1}}\sum_{j_1,\ldots,j_{n+1}=1}^d
\sum_{\scriptsize\begin{array}{c}k_1+\cdots+k_{n+1}=\alpha\\\ell_1+\cdots+\ell_{n+1}=\alpha\\
\sum_{s=1}^{n+1}(s-1)(\ell_s-k_s)=m\end{array}}
\prod_{a=1}^{n+1}(x_{j_a}^{k_a}\overline{x_{j_a}}^{\ell_a})
\binom{\alpha}{k_1,\ldots,k_{n+1}}\binom{\alpha}{\ell_1,\ldots,\ell_{n+1}}.$$
\end{theorem}

\subsection{Exact formulas}
The second of these directions to obtain an exact formula for
$\mu_T^{2\alpha}(n;m)$ for given $\alpha$ and $m$. By definitions,
it is clear that $\alpha!^2\mu_T^{2\alpha}(n;m)$ is the
coefficient of $w^m$ in $E_T^{\alpha,\alpha}(x,w)$ (see
Section~\ref{sec22} for definitions). Therefore, similarly as in
Section~\ref{sec22}, if applying $\mu$-algorithm with finding the
coefficient of $w^m$ (instead the free coefficient of $w$, see the
third step of the algorithm), then we get the following result.

\begin{corollary}
Let $T=\{x_1,x_2,\ldots,x_d\}$ be any finite set of complex
numbers. Then
$$\mu_T^{2\alpha}(n;m)=\left\{
        \begin{array}{ll}
        \frac{n+1-|m|}{d^2}A_T^{0,1}A_T^{1,0} & m\neq0,\, -n\leq m\leq m\\
        \frac{n+1}{d}A_T^{1,1} & m=0\\
        $0$                     & \mbox{otherwise}.
        \end{array}
        \right.$$
\end{corollary}
This result can be extended to the case of $\mu_T^4(n)$. Since the
answers become very cumbersome, we present here only the simplest
case when $A_T^{1,0}=\sum_{j=1}^dx_j=0$.

\begin{corollary}
Let $T=\{x_1,x_2,\ldots,x_d\}$ be any finite set of complex
numbers such that $A_T^{1,0}=\sum_{j=1}^dx_j=0$. Then
$$\mu_T^{2\alpha}(n;m)=\left\{
        \begin{array}{ll}
        \frac{n+1-|m/2|}{d^2}A_T^{0,2}A_T^{2,0} & m\neq0,\, -2n\leq m\leq2n,\, m\mbox{ even}\\
        \frac{n+1}{d}A_T^{1,1}+\frac{4}{d^2}\binom{n+1}{2}(A_T^{1,1})^2 & m=0\\
        $0$                     & \mbox{otherwise}.
        \end{array}
        \right.$$
\end{corollary}

\subsection{Average integrals}
The third of these directions is to consider the general case to
find an explicit formula for $e_T(n,s,t,m)$. Theorem~\ref{maa} and $\mu$-algorithm
gives a complete answer for the generating function for these
numbers, and explicit formula for $e_T(n,s,t,m)$. For example, the
following result is true.

\begin{corollary}
Let $T=\{x_1,x_2,\ldots,x_d\}$ be any finite set of complex
numbers. Then for any $n\geq0$ and $m\in\mathbb{Z}$,
$$e_T(n,1,2,m)=\left\{
               \begin{array}{ll}
               \frac{1}{d}A_T^{1,2}& 0\leq m\leq n\\
               0                    & \mbox{otherwise}
               \end{array}
               \right.\qquad\mbox{and}\qquad
e_T(n,2,1,m)=\left\{
               \begin{array}{ll}
               \frac{1}{d}A_T^{2,1}& -n\leq m\leq 0\\
               0                    & \mbox{otherwise}
               \end{array}
               \right.
.$$
\end{corollary}


\end{document}